\documentclass[11pt]{amsart}
\usepackage{amsthm}
\usepackage{amssymb}
\usepackage{amsmath}
\newtheorem{lemma}{Lemma}[section]
\newtheorem{theorem}[lemma]{Theorem}

\newtheorem{proposition}{Proposition}[section]

\newtheorem{corollary}{Corollary}[section]
\newtheorem*{corollary*}{Corollary}

\numberwithin{equation}{section}

\gdef\myletter{}

\let\savetheequation\theequation
\def\theequation{\savetheequation\myletter}

\newcommand{\CC}{{\mathbb C}}

\newcommand{\RR}{{\mathbb R}}

\def \bar{\overline}

\begin{document}

\vskip 3mm

\title[Christoffel function asymptotics ]{\bf Strong asymptotics for Christoffel functions of planar measues }  

\author{T. Bloom* and N. Levenberg}{\thanks{*Supported in part by an NSERC of Canada grant}}

\address{University of Toronto, Toronto, Ontario M5S 2E4 Canada}  
\email{bloom@math.toronto.edu}

\address{Indiana University, Bloomington, IN 47405 USA}

\email{nlevenbe@indiana.edu}

\date{September 12, 2007}

\begin{abstract}
We prove a version of strong asymptotics of Christoffel functions with varying weights for a
general class of sets $E$ and measures $\mu$ in the complex plane $\CC$. This class includes all
regular measures $\mu$ in the sense of Stahl-Totik \cite{stahltotik} on regular compact sets $E$ in $\CC$ and
even allows varying weights. Our main theorems cover some known results for $E \subset \RR$, the
real line; in particular, we recover information in the case of $E=\RR$ with Lebesgue measure $dx$ and weight $w(x) = \exp (-Q(x))$ where $Q(x)$ is a nonnegative,
even degree polynomial having positive leading coefficient.\end{abstract}

\maketitle

\section{\bf Introduction.}
	\label{sec:intro}
	Let $w(x)$ be a positive uppersemicontinuous (usc) function on the interval $[-1,1]$. Let $\{q_j\}_{j=1,2,...}$ be a sequence of orthonormal polynomials with respect to the measure $d\nu(x)=w(x)dx$ where deg$q_j=j-1$. Form the sequence of Christoffel functions $K_n(z):=\sum_{j=1}^{n+1} |q_j(z)|^2$. It is straightforward to see that 
	$$ \frac{1}{2n}\log K_n(z) \to \log|z+\sqrt {z^2-1}| $$
	uniformly on $\CC$ as $n\to \infty$ and hence
	\begin{equation} \label{arcsine}d\mu_n(x):=[\Delta (\frac{1}{2n}\log K_n(x))]dx  \to \frac{1}{\pi \sqrt {1-x^2}}dx \end{equation} weak-* where $\Delta$ is the Laplacian. More generally, let $E$ be a compact subset of $\CC$, $w$ an admissible weight function on $E$, and $\mu$ a positive Borel measure on $E$ such that the triple $(E,w,\mu)$ satisfies a weighted Bernstein-Markov inequality (see (\ref{wtdbernmark})). If we take, for each $n=1,2,...$, a set of orthonormal polynomials $q_1^{(n)},...,q_{n+1}^{(n)}$ with respect to the varying measures $w(z)^{2n}d\mu(z)$ where deg$q_j^{(n)}=j-1$ and form the sequence of Christoffel functions $K_n(z):=\sum_{j=1}^{n+1} |q_j^{(n)}(z)|^2$, then the functions $\frac{1}{2n}\log K_n(z)$ converge uniformly on $\CC$ and 
	\begin{equation} \label{weight} d\mu_n:=\Delta (\frac{1}{2n}\log K_n)  \to d\mu_{eq}^w \end{equation} 
	weak-* where $\mu_{eq}^w$ is the potential-theoretic weighted equilibrium measure (cf., \cite{bloomnew}, Lemma 2.3).

	A deeper result than (\ref{arcsine}) for the interval $[-1,1]$ is a stronger asymptotic for $K_n$: 
	\begin{equation} \label{asymp} \lim_{n\to \infty} \frac{K_n(x)}{n+1} = \pi \sqrt {1-x^2}/w(x) \end{equation}
	a.e. on $(-1,1)$. In \cite{totik2}, Totik generalized (\ref{asymp}) to certain ``regular'' measures $d\nu(x)$ on the real line having compact support $E$ where $\CC\setminus E$ is regular for the Dirichlet problem (he later observed that the regularity assumption on $E$ was unnecessary; cf., \cite{totiknew}, section 8). Here, the arcsine measure $ \frac{1}{\pi \sqrt {1-x^2}}dx$ is replaced by the (unweighted) equilibrium measure $d\mu_{eq}(x)$ for $E$ and the right-hand-side of (\ref{asymp}) is replaced by the appropriate Radon-Nikodym derivative. An earlier result of Totik (in \cite{totik1}) gives an analogous generalization of (\ref{weight}) in the special case where $E$ is a finite union of intervals in the real line, $\mu$ is normalized Lebesgue measure $dx$ on $E$, and $w$ is a positive continuous function on $E$. That is, forming the sequence of Christoffel functions $K_n(x):=\sum_{j=1}^{n+1} |q_j^{(n)}(x)|^2$, in this case, the asymptotic relation
	\begin{equation} \label{thm22}\frac{1}{n+1}K_n(x)w(x)^{2n}dx  \to d\mu_{eq}^w(x) \ \hbox{weak-* }\end{equation} 
	holds. We will refer to any result similar to (\ref{thm22}) as {\it strong asymptotics of Christoffel functions with varying weights} (this is not standard terminology).
	
	Motivated by the study of statistical quantities related to distributions of eigenvalues of random matrices, Johansson \cite{joh} and others (cf., Pastur \cite{pastur} and related articles) studied strong asymptotics of Christoffel functions with varying weights in the situation where $E$ is the whole real line $\RR$ and $w(x)=\exp (-Q(x))$ where $Q(x)$ is an even degree polynomial with positive leading coefficient and $Q(x)\geq 0$ on $\RR$; e.g., $Q(x)=x^{2m}, \ m=1,2,...$. We recommend Chapter 6 of Deift's elegant book  \cite{deift}; for the reader's convenience we include a statement (Corollary 2.1) and proof of a result in this case. Recently there has been a flurry of activity in {\it universality} limits involving orthogonal polynomials on subsets of the real line -- roughly speaking, the same (strong) asymptotics occur for a wide class of measures supported on the same set $E\subset \RR$. Very precise limiting behavior involving the reproducing kernels
	$$K_n(\mu;x,y):=\sum_{j=1}^{n+1}q_j(x)q_j(y)$$
where $\{q_j\}_{j=1,2,...}$ is a sequence of orthonormal polynomials (with positive leading coefficients) with respect to a measure $d\mu(x)$ on $E\subset \RR$ have been studied for certain $E$ and $\mu$ by Lubinsky, Simon, Totik and others (cf., \cite{lubinsky}, \cite{totiknew} and \cite{simon}). We utilize very different techniques, motivated from several complex variables, to prove a general version of (\ref{thm22}) for {\it nonpolar  compact sets $E$ in the complex plane} $\CC$ and triples $(E,w,\mu)$ satisfying a weighted Bernstein-Markov inequality (Theorem 2.2). Our main tool is the correspondence between weighted potential theory in $\CC$ and pluripotential theory (the study of plurisubharmonic functions) in $\CC^2$ as developed in \cite{bloom}. We also make use of some generalizations of certain statistical quantities described in \cite{deift}.  A crucial difference between our approach and that of Johansson is that we use potential theoretic consequences of measures satisfying a Bernstein-Markov inequality to prove the existence of ``free energy'' (Theorem 2.1). Our ``large deviation'' estimate (Proposition 4.2) is then a straightforward consequence of this result and the existence of the limit in (2.1); this latter item is a standard fact in potential theory. An alternate approach to studying strong asymptotics of Christoffel functions with varying weights in $\CC^N, \ N>1$ has been developed by Berman (cf., \cite{berman1}, \cite{berman2}). 

To keep the article self-contained, we include some brief background material on weighted potential theory in $\CC$ and pluripotential theory in $\CC^2$. We refer the reader to \cite{bloom} for details of stated results. For more on general univariate weighted potential theory, we refer the reader to \cite{safftotik}; for more on general pluripotential theory, we refer the reader to \cite{klimek}.  

In the next section, we give the relevant definitions and state our main results, Theorems 2.1 and 2.2. We prove Theorem 2.1 in section 3 and Theorem 2.2 in section 4. Corollary 2.1 is proved in section 5.

We would like to thank Vilmos Totik for his valuable comments; in particular, for kindly pointing out several recent papers on universality. 

\vskip 5mm

\section{\bf Main results.}
	\label{sec:KT}

		In this paper, we let $E\subset \CC$ be a closed set with an admissible weight function $w$: $w$ is a nonnegative, uppersemicontinuous function on $E$ with $\{z\in E:w(z)>0\}$ nonpolar (in particular, $E$ is nonpolar); if $E$ is unbounded we require also the growth condition $|z|w(z)\to 0$ as $|z|\to \infty$. It turns out that $S_w$ is always bounded (Theorem 1.3, p. 27 of \cite{safftotik}). 
				
		The limit 
	\begin{equation} \label{wtdvdm} \lim_{n\to \infty} \bigl[\max_{\lambda_i\in E}|VDM(\lambda_0,...,\lambda_n)|w(\lambda_0)^n\cdots w(\lambda_n)^n\bigr] ^{2/n^2}:=\delta^w(E) \end{equation} 
	exists and is called the weighted transfinite diameter of $E$ (with respect to $w$). Here $VDM(\zeta_1,...,\zeta_n)=\det [\zeta_i^{j-1}]_{i,j=1,...,n}=\prod_{j<k}(\zeta_j-\zeta_k)$ is the classical Vandermonde determinant. Points $\lambda_0,...,\lambda_n\in E$ for which $$|VDM(\lambda_0,...,\lambda_n)|w(\lambda_0)^n\cdots w(\lambda_n)^n$$
$$= |\det
\left[
\begin{array}{ccccc}
 1 &\lambda_0 &\ldots  &\lambda_0^n\\
  \vdots  & \vdots & \ddots  & \vdots \\
1 &\lambda_n &\ldots  &\lambda_n^n
\end{array}
\right]| \cdot w(\lambda_0)^n\cdots w(\lambda_n)^n$$
 is maximal are called {\it weighted Fekete points of order $n$}. The quantity $\delta^w(E)$ in (\ref{wtdvdm}) comes from a discrete version of a weighted energy minimization problem: for a probability measure $\tau $ on $E$, consider the weighted energy
	$$I^w(\tau):=\int_E \int_E \log \frac {1}{|z-t|w(z)w(t)}d\tau(t)d\tau(z).$$
	Then 
	\begin{equation} \label{wtdenmin} \inf_{\tau}I^w(\tau)=-\log \delta^w(E)\end{equation} 
	where the infimum is taken over all probability measures $\tau $ on $E$. Moreover, the infimum is attained by a unique measure $\mu_{eq}^w$. If $w\equiv 1$ we are in the classical (unweighted) case and we simply write $\mu_{eq}$. We remark that there exists $\eta >0$ such that the support $S_w$ of $\mu_{eq}^w$ is contained in $\{z\in E: w(z) \geq \eta\}$ (Remark 1.4, p. 27 of \cite{safftotik}).
			
		A weighted polynomial on $E$ is a function of the form $w(z)^np_n(z)$ where $p_n$ is a holomorphic polynomial of degree at most $n$. Let $\mu$ be a measure with support in $E$ such that $(E,w,\mu)$ satisfies a Bernstein-Markov inequality for weighted polynomials (referred to as a {\it weighted B-M inequality} in \cite{bloom}): given $\epsilon >0$, there exists a constant $M=M(\epsilon)$ such that for all weighted polynomials $w^np_n$
	\begin{equation} \label{wtdbernmark}||w^np_n||_E\leq M(1+\epsilon)^n||w^np_n||_{L^2(\mu)}.\end{equation}
	In this setting, we will restrict our attention to compact sets $E$. 
	
	For a compact set $K\subset \CC^N$ and a measure $\nu$ on $K$, we say that the pair $(K,\nu)$ satisfies the Bernstein-Markov inequality for holomorphic polynomials in $\CC^N$ if, given $\epsilon >0$, there exists a constant $\tilde M=\tilde M(\epsilon)$ such that for all such polynomials $Q_n$
	\begin{equation} \label{bernmark} ||Q_n||_K\leq \tilde M(1+\epsilon)^n||Q_n||_{L^2(\nu)}. \end{equation}
	The terminology ``Bernstein-Markov'' in this context is standard in several complex variables. Our main results are as follows.
	 
\begin{theorem} \label{zntheorem} Let $E$ be compact and let $(E,w,\mu)$ satisfy a Bernstein-Markov inequality for weighted polynomials. Then
$$\lim_{n\to \infty} Z_n^{1/n^2} = \delta^w(E)$$ where
\begin{equation} \label{zn} Z_n=Z_n(E,w,\mu):= \end{equation}
$$\int_{E^{n+1}} |VDM(\lambda_0,...,\lambda_n)|^2w(\lambda_0)^{2n} \cdots w(\lambda_n)^{2n}d\mu(\lambda_0) \cdots d\mu(\lambda_n).$$
\end{theorem}

\begin{theorem} \label{rntheorem} Let $E$ be compact and let $(E,w,\mu)$ satisfy a Bernstein-Markov inequality for weighted polynomials. Define the probability measures $$d\mu_n(z):=\frac {1}{Z_n} R_1^{(n)}(z)w(z)^{2n}d\mu(z)$$ where 
\begin{equation} \label{r1} R_1^{(n)}(z):=\end{equation}
$$\int_{E^{n}} |VDM(\lambda_0,...,\lambda_{n-1},z)|^2w(\lambda_0)^{2n} \cdots w(\lambda_{n-1})^{2n}d\mu(\lambda_0) \cdots d\mu(\lambda_{n-1}).$$
Then $d\mu_n(z) \to d\mu^w_{eq}(z)$ weak-*.
\end{theorem}

Utilizing similar arguments as in the proof of Theorem \ref{zntheorem}, we obtain another proof (cf., \cite{deift}, Chapter 6) of the following.

\begin{corollary} \label{zndeift} Let $w(x)=\exp (-Q(x))$ where $Q(x)$ is a nonnegative, even degree polynomial on the real line $\RR$ having positive leading coefficient. Then
\begin{equation} \label{mathcalzn} \lim_{n\to \infty} {\mathcal Z}_n^{1/n^2} = \delta^w(\RR)\end{equation} where
$$ {\mathcal Z}_n={\mathcal Z}_n(\RR,w,dx):= $$
$$\int_{\RR^{n+1}} |VDM(\lambda_0,...,\lambda_n)|^2w(\lambda_0)^{2n} \cdots w(\lambda_n)^{2n}d\lambda_0 \cdots d \lambda_n.$$
\end{corollary}

\noindent {\bf Remark 2.1}. We observe that with the notation in (\ref{r1}) and (\ref{zn})
\begin{equation}
\label{r1form} \frac {R_1^{(n)}(z)}{Z_n} = \frac{1}{n+1} \sum_{j=1}^{n+1} |q_j^{(n)}(z)|^2 \end{equation}
where $q^{(n)}_1,...,q^{(n)}_{n+1}$ are orthonormal polynomials of degree $0,1,...,n$ with respect to the measure $w(z)^{2n}d\mu(z)$; hence Theorem \ref{rntheorem} generalizes (\ref{thm22}). To verify (\ref{r1form}), if we apply Gram-Schmidt in $L^2(w^{2n}\mu )$ to the monomials $1,z,...,z^{n}$ to obtain orthogonal polynomials $p^{(n)}_1(z)\equiv 1,...,p^{(n)}_{n+1}(z)$, we have, upon using elementary row operations on $VDM(\lambda_0,...,\lambda_{n-1},z)$ and expanding the integrand in (\ref{r1}), $$R_1^{(n)}(z)=$$
$$\sum_{I,S}\sigma(I)\cdot \sigma(S)\int_{E^{n}} p^{(n)}_{i_0}(\lambda_0)\cdots p^{(n)}_{i_n}(z)\overline {p^{(n)}_{s_0}(\lambda_0)}\cdots \overline {p^{(n)}_{s_n}(z)}w(\lambda_0)^{2n} \cdots $$
$$\cdots w(\lambda_{n-1})^{2n}d\mu(\lambda_0) \cdots d\mu(\lambda_{n-1})$$
$$=\sum_{I,S}\sigma(I)\cdot \sigma(S) \bigl[\int_E p^{(n)}_{i_0}(\lambda_0)\overline {p^{(n)}_{s_0}(\lambda_0)}w(\lambda_0)^{2n} d\mu(\lambda_0)\cdots $$
$$\cdots \int_E p^{(n)}_{i_{n-1}}(\lambda_{n-1})\overline {p^{(n)}_{s_{n-1}}(\lambda_{n-1})}w(\lambda_{n-1})^{2n}d\mu(\lambda_{n-1})\bigr] p^{(n)}_{i_n}(z) \overline {p^{(n)}_{s_n}(z)} $$
\begin{equation}
\label{r1top}=n!\sum_{j=1}^{n+1}(\prod_{i\not= j}||p^{(n)}_i||_{L^2(w^{2n}\mu )}^2)|p^{(n)}_j(z)|^2.\end{equation} 
Here $I=(i_0,...,i_n)$ and $S=(s_0,...,s_n)$ are permutations of $(0,1,...,n)$ and $\sigma(I)$ is the sign of $I$ ($+1$ if $I$ is even; $-1$ if $I$ is odd). Then
$$Z_n=\int_E R_1^{(n)}(z)w(z)^{2n}d\mu(z)$$
$$=n!\sum_{j=1}^{n+1}(\prod_{i\not= j}||p^{(n)}_i||_{L^2(w^{2n}\mu )}^2)\int_E |p^{(n)}_j(z)|^2w(z)^{2n}d\mu(z)$$
\begin{equation}
\label{r1bottom}=n!\sum_{j=1}^{n+1}(\prod_{i=1}^{n+1}||p^{(n)}_i||_{L^2(w^{2n}\mu )}^2)=(n+1)!\prod_{i=1}^{n+1}||p^{(n)}_i||_{L^2(w^{2n}\mu )}^2.\end{equation}
Dividing (\ref{r1top}) by (\ref{r1bottom}) yields (\ref{r1form}) since $|q^{(n)}_j(z)|=|p^{(n)}_j(z)|/||p^{(n)}_j||_{L^2(w^{2n}\mu )}$.

In particular, then, if we take $E$ to be a finite union of intervals on the real line, $d\mu(x) =dx=$ Lebesgue measure on $E$, and $w$ positive and continuous, then $(E,w,dx)$ satisfies a Bernstein-Markov inequality for weighted polynomials (see Remark 2.2 below). Theorem \ref{rntheorem} gives  
$$[\frac{1}{n+1} \sum_{j=1}^{n+1} |q^{(n)}_j(x)|^2]\cdot w(x)^{2n}dx \to \rho(x)dx$$
weak-* where $\rho(x)dx = d\mu^w_{eq}(x)$. For such $E$ and $w$, under the additional hypotheses that $\rho$ be continuous on an interval $J$ in Int$E$ -- this occurs if $J\subset \ $Int$(S_w)$ and $w$ is $C^{1+\epsilon}$ on a neighborhood of $J$ -- Totik \cite{totik1} proved that 
$$[\frac{1}{n+1} \sum_{j=1}^{n+1} |q^{(n)}_j(x)|^2]\cdot w(x)^{2n} \to \rho(x)$$
uniformly on $J$. 

For more general subsets $E\subset \RR$ and measures $\mu$ such that $(E,\mu)$ satisfies a Bernstein-Markov inequality, if $w\equiv 1$, Theorem 2.2 implies   
$$[\frac{1}{n+1} \sum_{j=1}^{n+1} |q_j(x)|^2]d\mu(x) \to d\mu_{eq}(x).$$
In \cite{totik2}, Totik proves this result for a regular compact subset $E$ of $\RR$ ($\CC \setminus E$ is regular for the Dirichlet problem) with a ``regular'' measure $\mu$. According to Theorem 3.2.3 of \cite{stahltotik}, for a regular compact set $E$, regularity of $\mu$ is equivalent to $(E,\mu)$ satisfying a Bernstein-Markov inequality.
\vskip 3mm

\noindent {\bf Remark 2.2}. If $E$ is a regular compact set in $\CC$, then $(E,\mu_{eq})$ satisfies the Bernstein-Markov inequality. More generally, for such $E$, if $w$ is an admissible  continuous weight function, then the triple $(E,w,\mu_{eq}^w)$ satisfies a Bernstein-Markov inequality for weighted polynomials (cf., Corollary 3.1 \cite{bloom}). We mention that Theorem 4.2.3 of \cite{stahltotik} provides a sufficient condition for the pair $(E,\mu)$ to satisfy a Bernstein-Markov inequality when $E= \ $supp$\mu\subset \CC$ is a regular compact set. For instance, a finite union $E$ of intervals is regular, and $(E,dx)$ satisfies this sufficient condition where $dx$ is Lebesgue measure on $E$. If $w$ is a positive, continuous weight on $E$, appealing to Theorem 3.2.3 (vi) of \cite{stahltotik} with $g_n=w^{2n}$, it follows that $(E,w,dx)$ satisfies a Bernstein-Markov inequality for weighted polynomials.

\vskip 3mm

\noindent {\bf Remark 2.3}. In the setting of \cite{deift}, where $E=\RR$ and $w(x)=\exp (-Q(x))$ with $Q(x)$ an even degree polynomial having positive leading coefficient and $Q(x)\geq 0$ on $\RR$, Corollary \ref{zndeift} is referred to as the existence of the {\it free energy} (Corollary 6.90 in \cite{deift}). Note that $w$ is an admissible weight. 
\vskip 5mm

\section{\bf Proof of Theorem 2.1.}
	\label{sec:KT}
	
	We let $E\subset \CC$ be a nonpolar compact set with an admissible weight function $w$. Following \cite{bloom}, we define the circled set 
	$$F=F(E,w):=\{(t,z)=(t,\lambda t)\in \CC^2: \lambda \in E, \ |t|=w(\lambda)\}.$$
	We first relate weighted univariate Vandermonde determinants for $E$ with homogeneous bivariate Vandermonde determinants for $F$. To this end, for each positive integer $n$, choose $n+1$ points $\{(t_i,z_i)\}_{i=0,...,n}$ in $F$ and form the $n-$homogeneous Vandermonde determinant
	$$VDMH_n((t_0,z_0),...,(t_{n},z_{n})):=\det \bigl[t_i^{n-j}z_i^j\bigr]_{i,j=0,...,n}.$$
	Note that we evaluate the $n+1$ homogeneous monomials 
	$$t^n,t^{n-1}z,...,tz^{n-1},z^n$$ 
	at the $n+1$ points $\{(t_i,z_i)\}_{i=0,...,n}$. Factoring $t_i^n$ out of the $i-$th row, we obtain
	$$VDMH_n((t_0,z_0),...,(t_{n},z_{n}))=t_0^n\cdots t_n^n\cdot VDM(\lambda_0,...,\lambda_n);$$
i.e.,
\begin{align} \label{matrix}
\left|
\begin{array}{ccccc}
 t_0^n &t_0^{n-1}z_0 &\ldots  &z_0^n\\
  \vdots  & \vdots & \ddots  & \vdots \\
t_n^n &t_n^{n-1}z_n &\ldots  &z_n^n
\end{array}
\right| 
= t_0^n\cdots t_n^n
\left|
\begin{array}{ccccc}
 1 &\lambda_0 &\ldots  &\lambda_0^n\\
  \vdots  & \vdots & \ddots  & \vdots \\
1 &\lambda_n &\ldots  &\lambda_n^n
\end{array}
\right|, 
\end{align}
	where $\lambda_j = z_j/t_j$ provided $t_j\not = 0$ and $VDM(\lambda_0,...,\lambda_n)=\prod_{j<k}(\lambda_j -\lambda_k)$ is a standard (univariate) Vandermonde determinant. By definition of $F$, since $(t_i,z_i)\in F$, we have $ |t_i|=w(\lambda_i)$ so that from (\ref{matrix})
	$$|VDMH_n((t_0,z_0),...,(t_{n},z_{n}))|=|VDM(\lambda_0,...,\lambda_n)|w(\lambda_0)^n\cdots w(\lambda_n)^n.$$
	Thus
	$$\max_{(t_i,z_i)\in F}|VDMH_n((t_0,z_0),...,(t_{n},z_{n}))|=$$
	$$\max_{\lambda_i\in E}|VDM(\lambda_0,...,\lambda_n)|w(\lambda_0)^n\cdots w(\lambda_n)^n.$$
	Note that the maximum will occur when all $|t_j|=w(\lambda_j)>0$. Now the limit
	$$\lim_{n\to \infty} \bigl[\max_{(t_i,z_i)\in F}|VDMH_n((t_0,z_0),...,(t_{n},z_{n}))|\bigr]^{1/n^2}=:D^H(F)$$
	exists (the {\it homogeneous bivariate transfinite diameter} $D^H(F)$ of $F$; cf., \cite{jed}); also, as previously mentioned, the limit 
	$$\lim_{n\to \infty} \bigl[\max_{\lambda_i\in E}|VDM(\lambda_0,...,\lambda_n)|w(\lambda_0)^n\cdots w(\lambda_n)^n\bigr] ^{2/n^2}:=\delta^w(E)$$
	exists (the weighted transfinite diameter of $E$ with respect to $w$) and thus we have 
	\begin{equation} \label{homvswtd} \delta^w(E)=D^H(F)^2. \end{equation}

	Let $\mu$ be a measure with support in $E$ such that $(E,w,\mu)$ satisfies a Bernstein-Markov inequality for weighted polynomials. Note that the integrand 
	$$|VDM(\lambda_0,...,\lambda_n)|^2w(\lambda_0)^{2n} \cdots w(\lambda_n)^{2n}$$ in the definition of $Z_n$ in (\ref{zn}) thus has a maximal value on $E^{n+1}$ whose $1/n^2$ root tends to $\delta^w(E)$. To show that the integrals themselves have the same property, we proceed as follows. On the set $F\subset \CC^2$, there exists a measure $\nu$ associated to $\mu$ such that $(F,\nu)$ satisfies the Bernstein-Markov property for holomorphic polynomials in $\CC^2$; i.e., (\ref{bernmark}) holds. Indeed, take
	$$\nu:= m_{\lambda} \otimes \mu, \ \lambda \in E$$
	where $m_{\lambda}$ is normalized Lebesgue measure on the circle $|t|=w(\lambda)$ in the complex $t-$plane given by  
	$$C_{\lambda}:=\{(t,t\lambda)\in \CC^2: t\in \CC\}.$$
	That is, if $\phi$ is continuous on $F$, 
	$$\int_F \phi(t,z) d\nu (t,z) = \int_E\bigl[\int_{C_{\lambda}}\phi(t,t\lambda)dm_{\lambda}(t)\bigr] d\mu(\lambda).$$
	Equivalently, if $\pi:\CC^2 \to \CC$ via $\pi(t,z)=z/t:=\lambda$, then $\pi_*(\nu) =\mu$. Moreover, if $p_1(t,z)$ and $p_2(t,z)$ are two homogeneous polynomials in $\CC^2$ of degree $n$, say, and we write $$p_j(t,z)=p_j(t,\lambda t)=t^np_j(1,\lambda)=:t^n G_j(\lambda), \ j=1,2$$
	for univariate $G_j$, then it is straightforward to see that 
	\begin{equation}
\label{orthog}	\int_F p_1(t,z)\overline {p_2(t,z)}d\nu(t,z) = \int_E G_1(\lambda)\overline {G_2(\lambda)}w(\lambda)^{2n}d\mu(\lambda) \end{equation}
	(cf., \cite{bloom}, Lemma 3.1 and its proof). Note that if $p(t,z)=t^iz^{n-i}$ for $i=0,...,n$, then 
	$$p(t,z)=t^n\bigl(z/t\bigr)^{n-i}=t^n G(\lambda)$$
	where $G(\lambda)= \lambda^{n-i}$.

	\begin{proposition} Let 
	$$\tilde Z_n:= \int_{F^{n+1}} |VDMH_n((t_0,z_0),...,(t_{n},z_{n}))|^2 d\nu(t_0,z_0) \cdots d\nu(t_n,z_n).$$
	Then $\tilde Z_n= Z_n$. 
	\end{proposition}
	
	\begin{proof} Using the notation from the proof of (\ref{r1form}),  expanding the homogeneous Vandermonde determinant in $\tilde Z_n$ gives
	$$\tilde Z_n =\sum_{I,S}\sigma(I)\cdot \sigma(S) \bigl[\int_F t_0^{i_0}z_0^{n-i_0}\bar t_0^{s_0} \bar z_0^{n-s_0}d\nu(t_0,z_0)\cdots $$
$$\cdots \int_F t_n^{i_n}z_n^{n-i_n}\bar t_n^{s_n} \bar z_n^{n-s_n}d\nu(t_n,z_n)\bigr].$$
Expanding the ordinary  Vandermonde determinant in $Z_n$ gives 
$$Z_n = \sum_{I,S}\sigma(I)\cdot \sigma(S) \bigl[\int_E\lambda_0^{n-i_0} \bar \lambda_0^{n-s_0}w(\lambda_0)^{2n}d\mu(\lambda_0)\cdots $$
$$\cdots \int_E\lambda_n^{n-i_n} \bar \lambda_n^{n-s_n}w(\lambda_n)^{2n}d\mu(\lambda_n)\bigr].$$
Since $|t_j|=w(\lambda_j)$, using (\ref{orthog}) completes the proof.
	\end{proof}
	
	We need the following result, Theorem 5.9 of \cite{jed}, which is the homogeneous analogue of Theorem 3.3 in \cite{BBCL}.
		
\begin{proposition}  \label{bbclhom} Let $F\subset \CC^2$ be a circled set and let $\nu$ be a measure with support in $F$ such that $(F,\nu)$ satisfies the Bernstein-Markov property. Then $$D^H(F)=\lim_{n\to \infty} {\mathcal G}_n^{1/2n^2}$$ where 
$${\mathcal G}_n:= \det \bigl[\int_F t^{n-j}z^j\bar t^{n-i}\bar z^i d\nu(t,z)\bigr]_{i,j=0,...,n}$$
is the $n-$th homogeneous Gram determinant associated with $(F,\nu)$. 

\end{proposition}
	
	The last step is to work in $\CC^2$ with the $\tilde Z_n$ integrals and verify the following.
	
	\begin{proposition} We have 
	$$\lim_{n\to \infty} \tilde Z_n^{1/2n^2}= D^H(F).$$	\end{proposition}
	
	\begin{proof} Fix $n$ and consider the monomials 
	$$t^n,t^{n-1}z,...,tz^{n-1},z^n$$
	utilized in $VDMH_n((t_0,z_0),...,(t_{n},z_{n}))$ and in the computation of the Gram determinant ${\mathcal G}_n$ associated to $(F,\nu)$. Use Gram-Schmidt in $L^2(\nu)$ to obtain orthogonal polynomials 
	$$p_0(t,z)=t^n, \ p_1(t,z) = t^{n-1}z + \cdots, ..., p_{n}(t,z) = z^{n} + \cdots.$$
	Then 
	$$VDMH_n((t_0,z_0),...,(t_{n},z_{n}))=\det \bigl[ p_i(t_j,z_j)\bigr]_{i,j=0,...,n}.$$
	By orthogonality, as in the calculation of $Z_n$ in proving (\ref{r1bottom}) in the introduction, we obtain
	$$Z_n =(n+1)! ||p_0||_{L^2(\nu)}^2\cdots ||p_n||_{L^2(\nu)}^2.$$
	On the other hand, using the orthogonal polynomials diagonalizes the $n-$th Gram matrix while preserving its determinant; hence 
	$${\mathcal G}_n=||p_0||_{L^2(\nu)}^2\cdots ||p_n||_{L^2(\nu)}^2.$$ 
	From Proposition \ref{bbclhom}, we have
	$$\lim_{n\to \infty}(||p_0||_{L^2(\nu)}^2\cdots ||p_n||_{L^2(\nu)}^2)^{1/2n^2} =D^H(F)$$
	and the result follows.
	\end{proof}
	
	Combining Propositions 3.1 and 3.3 with equation (\ref{homvswtd}) completes the proof of Theorem \ref{zntheorem}. \hfill $\Box$

\smallskip

	\noindent {\bf Remark 3.1}. A version of Theorem \ref{zntheorem} is valid in $\CC^N$ for $N>1$. Here, in the definition of $Z_n$ in (\ref{zn}), we replace $VDM(\lambda_1,...,\lambda_n)$ by the generalized Vandermonde determinant
$$VDM(\lambda_1,...,\lambda_{h_n}):=|\det [e_i(\lambda_j)]_{i,j=1,...,h_n}|$$
where $e_1(z),...,e_i(z),...,e_{h_n}(z)$ is a listing of the monomials in
${\bf C}^N$ of degree at most $n$. If $\ell_n= \sum_{i=1}^{h_n}$deg$e_i$, the result is then that
$$\lim_{n\to \infty} Z_N^{1/2\ell_n}=\delta^w(E) $$
where in (\ref{wtdvdm}) we use
$$\delta^w(E) := \lim_{n\to \infty} \bigl[\max_{\lambda_i\in E}|VDM(\lambda_1,...,\lambda_{h_n})|w(\lambda_1)^n\cdots w(\lambda_{h_n})^n\bigr] ^{1/\ell_n}.$$
Details of this and related results will be given in a future work.

\vskip 5mm

\section{\bf Proof of Theorem 2.2.}
	\label{sec:KT}
	
	As in the previous section, we let $E\subset \CC$ be a nonpolar compact set with an admissible weight function $w$. Let $\mu$ be a measure with support in $E$ such that $(E,w,\mu)$ satisfies a Bernstein-Markov inequality for weighted polynomials. We will need the following fact, which is claimed in Remark 1.4 on p. 147 of \cite{safftotik}. This says that for any doubly indexed array of points $\{z_k^{(n_j)}\}_{k=1,...,n_j; \ j=1,2,...}$ in $E$  
which satisfies asymptotically the relation (\ref{wtdvdm}), the limiting measures 
\begin{equation} \label{wtdasymp} d\mu_{n_j} :=\frac{1}{n_j} \sum_{k=1}^{n_j}\delta_{z_k^{(n_j)}} \end{equation} have the same weak-* limit, the weighted equilibrium measure $d\mu^w_{eq}$.
	
	\begin{proposition} \label{asympfek} Let $E\subset \CC$ be compact and let $w$ be an admissible weight on $E$. If, for a subsequence of positive integers $\{n_j\}$ with $n_j \uparrow \infty$, the points $z_1^{(n_j)},...,z_{n_j}^{(n_j)}\in E$ are chosen so that 
$$\lim_{j\to \infty} \bigl[|VDM(z_1^{(n_j)},...,z_{n_j}^{(n_j)})|^2w(z_1^{(n_j)})^{2n_j} \cdots w(z_{n_j}^{(n_j)})^{2n_j}\bigr]^{1/n_j^2}=\delta^w(E),$$
then $d\mu_{n_j}\to d\mu^w_{eq}$ weak-* where $d\mu_{n_j}$ is defined in (\ref{wtdasymp}).	
	\end{proposition}
	
	We recall that the support $S_w$ of $\mu_{eq}^w$ is contained in $E_{\eta}:=\{z\in E: w(z) \geq \eta\}$ for some $\eta >0$ and we mention that the {\it proof} on p. 146 of \cite{safftotik} for weighted Fekete points works verbatim under the 
assumption that the points $z_1^{(n_j)},...,z_{n_j}^{(n_j)}$ lie in $E_{\tilde \eta}$ for some $\tilde \eta >0$. For the reader's convenience, and since we do not make this assumption, we include a proof of Proposition 4.1.

\begin{proof} 	Take a subsequence of the measures $\{\mu_{n_j}\}$ which converges weak-* to a probability measure $\sigma$ on $E$. We use the same notation for the subsequence and the original sequence. We show that $I^w(\sigma) = -\log \delta^w$; by uniqueness of the weighted energy minimizing measure (\ref{wtdenmin}) we will then have $\sigma = \mu^w_{eq}$. First of all, choose continuous admissible weight functions $\{w_m\}$ with $w_m \downarrow w$ and $w_m \geq \alpha_m >0$ on $E$ and for a real number $M$ let 
	$$h_{M,m}(z,t):=\min [M, \log \frac {1}{|z-t|w_m(z)w_m(t)}]\leq  \log \frac {1}{|z-t|w_m(z)w_m(t)}$$
	and
	$$h_M(z,t):=\min [M, \log \frac {1}{|z-t|w(z)w(t)}]\leq \log \frac {1}{|z-t|w(z)w(t)}.$$
	Then $h_{M,m}\leq h_M$. By the Stone-Weierstrass theorem, every continuous function on $E\times E$ can be uniformly approximated by finite sums of the form $\sum_j f_j(z)g_j(t)$ where $f_j,g_j$ are continuous on $E$; hence $\mu_{n_j}\times \mu_{n_j}\to \sigma \times \sigma$ and we have
$$I^w(\sigma)=\lim_{M\to \infty}\lim_{m\to \infty}\int_E \int_E h_{M,m}(z,t)d\sigma(z)d\sigma(t)$$
$$=\lim_{M\to \infty}\lim_{m\to \infty}\lim_{j \to \infty}\int_E \int_E h_{M,m}(z,t)d\mu_{n_j}(z)d\mu_{n_j}(t)$$
$$\leq \lim_{M\to \infty}\limsup_{j \to \infty}\int_E \int_E h_{M}(z,t)d\mu_{n_j}(z)d\mu_{n_j}(t)$$
since $h_{M,m}\leq h_M$. 
Now  
$$h_{M}(z_k^{(n_j)},z_l^{(n_j)})\leq \log \frac {1}{|z_k^{(n_j)}-z_l^{(n_j)}|w(z_k^{(n_j)})w(z_l^{(n_j)})}$$
if $k\not = l$ and hence
$$\int_E \int_E h_{M}(z,t)d\mu_{n_j}(z)d\mu_{n_j}(t)\leq $$
$$\frac{1}{n_j}M +(\frac{1}{n_j^2-n_j})\bigl[\sum_{k\not = l}\log \frac {1}{|z_k^{(n_j)}-z_l^{(n_j)}|w(z_k^{(n_j)})w(z_l^{(n_j)})}\bigr].$$
By assumption, given $\epsilon >0$, 
$$(\frac{1}{n_j^2-n_j})\bigl[\sum_{k\not = l}\log \frac {1}{|z_k^{(n_j)}-z_l^{(n_j)}|w(z_k^{(n_j)})w(z_l^{(n_j)})}\bigr]\leq -\log [\delta^w(E)-\epsilon]$$
for $j\geq j(\epsilon)$; in particular, $w(z_k^{(n_j)})> 0$ for such $j$ and hence
$$I^w(\sigma)\leq \lim_{M\to \infty}\limsup_{j \to \infty} \frac{1}{n_j}M -\log [\delta^w(E)-\epsilon] =-\log [\delta^w(E)-\epsilon]$$
for all $\epsilon >0$; i.e., $I^w(\sigma)=-\log \delta^w(E)$. \end{proof}

	We also need a ``large deviation'' result, which follows easily from Theorem 2.1. Define a probability measure ${\mathcal P}_n$ on $E^{n+1}$ via, for a Borel set $A\subset E^{n+1}$,
	$${\mathcal P}_n(A):=\frac{1}{Z_n}\int_{A} |VDM(z_0,...,z_n)|^2w(z_0)^{2n} \cdots w(z_n)^{2n}d\mu(z_0) \cdots d\mu(z_n).$$
		
\begin{proposition} \label{largedev} Given $\eta >0$, define
$$A_{n,\eta}:=$$
$$\{(z_0,...,z_n)\in E^{n+1}: |VDM(z_0,...,z_n)|^2w(z_0)^{2n} \cdots w(z_n)^{2n} \geq (\delta^w(E) -\eta)^{n^2}\}.$$
Then there exists $n^*=n^*(\eta)$ such that for all $n>n^*$, 
$${\mathcal P}_n(E^{n+1}\setminus A_{n,\eta})\leq (1-\frac{\eta}{2\delta^w(E)})^{n^2}.$$
	\end{proposition}	
	
	\begin{proof} From Theorem 2.1, given $\epsilon >0$, 
	$$Z_n \geq [\delta^w(E) -\epsilon]^{n^2}$$
	for $n\geq n(\epsilon)$. Thus
	$${\mathcal P}_n(E^{n+1}\setminus A_{n,\eta})=$$
	$$\frac{1}{Z_n}\int_{E^{n+1}\setminus A_{n,\eta}} |VDM(z_0,...,z_n)|^2w(z_0)^{2n} \cdots w(z_n)^{2n}d\mu(z_0) \cdots d\mu(z_n)$$
	$$\leq \frac{  [\delta^w(E) -\eta]^{n^2}} { [\delta^w(E) -\epsilon]^{n^2}}$$
	if $n\geq n(\epsilon)$. Choosing $\epsilon < \eta/2$ and $n^*=n(\epsilon)$ gives the result.
	\end{proof}
	
	To prove Theorem 2.2, we fix $\phi\in C(E)$. Recalling that 
	$$d\mu_n(z):=\frac {1}{Z_n} R_1^{(n)}(z)w(z)^{2n}d\mu(z),$$
for each $n$ we have
	$$\int_E\phi(z)d\mu_n(z)$$
	$$= \frac{1}{Z_n}\int_E \phi(z) \bigl(\int_{E^{n}} |VDM(z_0,...,z_{n-1},z)|^2w(z_0)^{2n} \cdots w(z_{n-1})^{2n}$$
	$$d\mu(z_0) \cdots d\mu(z_{n-1})\bigr)w(z)^{2n}d\mu(z)$$
	$$= \frac{1}{Z_n}\int_{E^{n+1}} \phi(z_n)  |VDM(z_0,...,z_n)|^2w(z_0)^{2n} \cdots w(z_n)^{2n}d\mu(z_0) \cdots d\mu(z_{n})$$
	$$=\frac{1}{Z_n}\int_{E^{n+1}} \frac{\sum_{j=0}^n\phi(z_j)}{n+1} |VDM(z_0,...,z_n)|^2w(z_0)^{2n} \cdots w(z_n)^{2n}d\mu(z_0) \cdots d\mu(z_{n})$$
	$$=:\int_{E^{n+1}}\psi_n(z_0,...,z_n)d{\mathcal P}_n(z_0,...,z_n)$$
	where $\psi_n(z_0,...,z_n):= \frac{\sum_{j=0}^n\phi(z_j)}{n+1}$. 
	
	Now take a sequence $\{\eta_j\}$ with $\eta_j \downarrow 0$ and a corresponding sequence $\{n_j\}$ with $n_j \geq n^*(\eta_j)$ from Proposition \ref{largedev} and $n_j \uparrow \infty$. By choosing $n_j$ larger if necessary we may assume that
	\begin{equation} \label{cond} n_j^2\eta_j\uparrow \infty \ \hbox{so that} \  (1-\frac{\eta_j}{2\delta^w(E)})^{n_j^2}\to 0.\end{equation}
	Choose points $\tilde z_0^{(n_j)},...,\tilde z_{n_j}^{(n_j)}\in A_{n_j,\eta_j}$ with
	$$\psi_{n_j}(\tilde z_0^{(n_j)},...,\tilde z_{n_j}^{(n_j)})=\max_{(w_0,...,w_{n_j})\in A_{n_j,\eta_j} }\psi_{n_j}(w_0,...,w_{n_j}).$$
	If $|\phi|\leq M$ on $E$, then $|\psi_{n_j}|\leq M$ on $E^{n_j+1}$;  using the large deviation result, Proposition 4.2, and (\ref{cond}),
	$$\limsup_{j\to \infty} \int_E \phi(z)d\mu_{n_j}(z) =\limsup_{j\to \infty}\int_{E^{n_j+1}}\psi_{n_j}d{\mathcal P}_{n_j}$$
	$$=\limsup_{j\to \infty}[\int_{A_{n_j,\eta_j} }\psi_{n_j}d{\mathcal P}_{n_j}
+\int_{E^{n_j+1}\setminus A_{n_j,\eta_j} }\psi_{n_j}d{\mathcal P}_{n_j}]$$
$$\leq \limsup_{j\to \infty}  \bigl(\frac{1}{n_j+1}\sum_{k=0}^{n_j}\phi(\tilde z_k^{(n_j)})+M(1-\frac{\eta_j}{2\delta^w(E)})^{n_j^2}\bigr)$$
	$$= \limsup_{j\to \infty}  \frac{1}{n_j+1}\sum_{k=0}^{n_j}\phi(\tilde z_k^{(n_j)}).$$
	Now since $\tilde z_0^{(n_j)},...,\tilde z_{n_j}^{(n_j)}\in A_{n_j,\eta_j}$, 
	$$|VDM(\tilde z_0^{(n_j)},...,\tilde z_{n_j}^{(n_j)})|^2w(\tilde z_0^{(n_j)})^{2n_j} \cdots w(\tilde z_{n_j}^{(n_j)})^{2n_j}\geq (\delta^w(E) -\eta_j)^{n_j^2}$$
	so that
	$$\lim_{j\to \infty}\bigl(|VDM(\tilde z_0^{(n_j)},...,\tilde z_{n_j}^{(n_j)})|^2w(\tilde z_0^{(n_j)})^{2n_j} \cdots w(\tilde z_{n_j}^{(n_j)})^{2n_j}\bigr)^{1/n_j^2}=\delta^w(E).$$
	By Proposition \ref{asympfek}, 
	$$\frac{1}{n_j+1} \sum_{k=0}^{n_j}\delta_{\tilde z_k^{(n_j)}}\to d\mu^w_{eq}.$$
	Thus
	$$\frac{1}{n_j+1}\sum_{k=0}^{n_j}\phi(\tilde z_k^{(n_j)})\to \int_E \phi(z)d\mu^w_{eq}(z)$$ and hence 
	\begin{equation} \label{weak} 	
	\limsup_{j\to \infty} \int_E \phi(z)d\mu_{n_j}(z) \leq  \int_E \phi(z)d\mu^w_{eq}(z). \end{equation}
	Applying (\ref{weak}) to $-\phi$ we obtain
$$\limsup_{j\to \infty} \int_E (-\phi(z))d\mu_{n_j}(z) \leq  \int_E (-\phi(z))d\mu^w_{eq}(z);$$
i.e., 
$$\liminf_{j\to \infty} \int_E \phi(z)d\mu_{n_j}(z) \geq  \int_E \phi(z)d\mu^w_{eq}(z),$$
so that 
\begin{equation} \label{weak2} \lim_{j\to \infty} \int_E \phi(z)d\mu_{n_j}(z) =  \int_E \phi(z)d\mu^w_{eq}(z). \end{equation}
Thus for any sequence of positive integers increasing to infinity we can choose a subsequence $\{n_j\}$ satisfying (\ref{cond}) for some $\eta_j\downarrow 0$ so that (\ref{weak2}) holds; hence $d\mu_n(z) \to d\mu^w_{eq}(z)$ weak-*.
\hfill $\Box$
\smallskip
	
\noindent {\bf Remark 4.1}.  More generally, if we consider, for any positive integer $m\geq 1$, the generalized {\it $m-$point correlation functions} $R_m^{(n)}(z_1,...,z_m)$ defined as 
$$R_m^{(n)}(z_1,...,z_m):=\int_{E^{n-m+1}} |VDM(\lambda_0,...,\lambda_{n-m},z_1,...,z_m)|^2\cdot $$
$$w(\lambda_0)^{2n} \cdots w(\lambda_{n-m})^{2n}d\mu(\lambda_0) \cdots d\mu(\lambda_{n-m}),$$
then one may verify that 
$$\frac {1}{Z_n} R_m^{(n)}(z_1,...,z_m)w(z_1)^{2n}\cdots w(z_m)^{2n}d\mu(z_1)\cdots d\mu(z_m)$$
converge weak-* as $n\to \infty$ to $d\mu^w_{eq}(z_1) \cdots d\mu^w_{eq}(z_m)$. See \cite{deift} for the case $E=\RR$ and $w(x)=\exp (-Q(x))$ where $Q(x)$ is an even degree polynomial with positive leading coefficient and $Q(x)\geq 0$ on $\RR$. In our setting, one may prove the analogue of Lemma 6.77 of \cite{deift} with slight modifications and then the proof of the analogues of Corollary 6.94 and Theorem 6.96 follow word-for-word.

\vskip 5mm

\section{\bf Proof of Corollary 2.1.}
	\label{sec:KT}

We indicate the modifications needed to prove Corollary \ref{zndeift}. We have $E=\RR$, $\mu = dx$ and $w(x)=\exp (-Q(x))$ with $Q(x)$ an even degree polynomial having positive leading coefficient and $Q(x)\geq 0$ on $\RR$. As mentioned in section 2, it is known that for unbounded sets and admissible measures, the support $S_w$ of the weighted energy minimizing measure $\mu_{eq}^w$ is compact. 
Related to this is the observation that the $L^2-$norms of our weighted polynomials essentially ``live'' on a compact subset of $\RR$ (cf., Theorem III.6.1 of \cite{safftotik}). To be precise, we can take $\tilde E$ to be a large enough compact interval in $\RR$ so that $S_w\subset \tilde E$ and such that there exist positive constants $a$ and $b$ independent of $n$ and $p_n$ so that if $p_n$ is a polynomial of degree at most $n$,
\begin{equation} \label{l2normlive} \int_{\RR} |p_n(x)|^2 w(x)^{2n}dx \leq (1+ae^{-bn})\int_{\tilde E} |p_n(x)|^2 w(x)^{2n}dx. \end{equation}

We apply Gram-Schmidt in $L^2(w^{2n}dx )$ to the monomials $1,x,...,x^{n}$ to obtain orthogonal polynomials $p^{(n)}_1(x)\equiv 1,...,p^{(n)}_{n+1}(x)$; and we apply Gram-Schmidt in $L^2(w^{2n}dx|_{\tilde E} )$ to the monomials $1,x,...,x^{n}$ to obtain orthogonal polynomials $q^{(n)}_1(x)\equiv 1,...,q^{(n)}_{n+1}(x)$. From (\ref{l2normlive}), 
\begin{equation} \label{l2norms} 
||p^{(n)}_j||^2_{L^2(w^{2n}dx )} \leq (1+ae^{-b(j-1)})||p^{(n)}_j||^2_{L^2(w^{2n}dx|_{\tilde E} )} \ \hbox{and} \end{equation}
$$ ||q^{(n)}_j||^2_{L^2(w^{2n}dx )} \leq (1+ae^{-b(j-1)})||q^{(n)}_j||^2_{L^2(w^{2n}dx|_{\tilde E} )}. $$
Also, from the definitions of the orthogonal polynomials, we have
\begin{equation} \label{l2normcomp} 
||q^{(n)}_j||_{L^2(w^{2n}dx|_{\tilde E} )}\leq ||p^{(n)}_j||_{L^2(w^{2n}dx|_{\tilde E} )} \leq ||p^{(n)}_j||_{L^2(w^{2n}dx )} \leq ||q^{(n)}_j||_{L^2(w^{2n}dx)}. \end{equation}
so that, combining (\ref{l2norms}) and (\ref{l2normcomp}), 
\begin{equation} \label{finalcomp} 
||q^{(n)}_j||_{L^2(w^{2n}dx|_{\tilde E} )}\leq ||p^{(n)}_j||_{L^2(w^{2n}dx|_{\tilde E} )}\leq (1+ae^{-b(j-1)})^{1/2}||q^{(n)}_j||_{L^2(w^{2n}dx|_{\tilde E} )}. \end{equation}

As in Remark 2.1, 
$$ {\mathcal Z}_n={\mathcal Z}_n(\RR,w,dx):= $$
$$\int_{\RR^{n+1}} |VDM(\lambda_0,...,\lambda_n)|^2w(\lambda_0)^{2n} \cdots w(\lambda_n)^{2n}d \lambda_0 \cdots d \lambda_n$$
can be written as (see (\ref{r1bottom})) 
$${\mathcal Z}_n = (n+1)!\prod_{i=1}^{n+1}||p^{(n)}_i||_{L^2(w^{2n}dx )}^2.$$
Note the $L^2-$norms are finite because of the decay as $|x|\to \infty$ of $w(x)$. Using (\ref{l2norms}), 
$$||p^{(n)}_j||^2_{L^2(w^{2n}dx )} \leq (1+ae^{-b(j-1)})||p^{(n)}_j||^2_{L^2(w^{2n}dx|_{\tilde E} )}\leq (1+ae^{-b(j-1)})||p^{(n)}_j||^2_{L^2(w^{2n}dx )};$$
multiplying these inequalities for $j=1,...,n+1$ and taking $n^2-$roots, we see that
\begin{equation} \label{mathcalznmod}\lim_{n\to \infty} {\mathcal Z}_n^{1/n^2} = \lim_{n\to \infty}\bigl[(n+1)!\prod_{i=1}^{n+1}||p^{(n)}_i||_{L^2(w^{2n}dx|_{\tilde E} )}^2\bigr]^{1/n^2}\end{equation}
provided this limit exists.

On the other hand, applying Theorem \ref{zntheorem} and (\ref{r1bottom}) to $(\tilde E,w, dx|_{\tilde E})$ (recall Remark 2.2), we have
\begin{equation} \label{tildee} \lim_{n\to \infty}\bigl[(n+1)!\prod_{i=1}^{n+1}||q^{(n)}_i||_{L^2(w^{2n}dx|_{\tilde E} )}^2\bigr]^{1/n^2}=\delta^w(\tilde E). \end{equation}
But since $\mu_{eq}^w$ has support in $\tilde E$, from (\ref{wtdenmin}), 
\begin{equation} \label{dwe=dwr}\delta^w(\tilde E)=\delta^w(\RR). \end{equation}
 The proof of (\ref{mathcalzn}), including the existence of the limit, now follows from (\ref{mathcalznmod}), (\ref{tildee}) and (\ref{dwe=dwr}) by applying (\ref{finalcomp}).

\hfill $\Box$

\end{document}